\documentclass{article}

\newtheorem{fed}{\textbf{Definition}}[section]
\newtheorem{thm}[fed]{\textbf{Theorem}}
\newtheorem{lemma}[fed]{\textbf{Lemma}}

\newtheorem{prop}[fed]{\textbf{Proposition}}

\usepackage{amssymb,bbm,graphicx,epsfig,psfrag,epic,eepic,latexsym}
\usepackage{amsmath}
\usepackage{mathrsfs}
\usepackage[dvips]{color}

\begin{document}
\title{The Gradient flow equation of Rabinowitz action functional in a symplectization}
\author{Urs Frauenfelder} \maketitle

\begin{abstract}
Rabinowitz action functional is the Lagrange multiplier functional of the negative area
functional to a constraint given by the mean value of a Hamiltonian. In this note we
show that on a symplectization there is a one-to-one correspondence between gradient
flow lines of Rabinowitz action functional and gradient flow lines of the restriction
of the negative area functional to the constraint. In the appendix we explain the motivation
behind this result. Namely that the restricted functional satisfies Chas-Sullivan additivity
for concatenation of loops which the Rabinowitz action functional does in general not do. 
\end{abstract}

\section{Introduction}

Assume that $M$ is a manifold, and $f,h \colon M \to \mathbb{R}$ are two smooth functions
such that $0$ is a regular value of $h$. Then critical points of $f$ restricted to
the hypersurface $h^{-1}(0) \subset M$ can be alternatively detected by the Lagrange 
multiplier functional
$$F \colon M \times \mathbb{R} \to \mathbb{R}, \quad (x,\tau) \mapsto f(x)+\tau h(x).$$
Indeed, if $(x,\tau)$ is a critical point of $F$, then $x$ is a critical point of
$f|_{h^{-1}(0)}$ and $\tau$ is referred to as the Lagrange multiplier. In \cite{frauenfelder}
it is explained how the Morse homology of the Lagrange multiplier functional coincides
with the singular homology of the constraint $h^{-1}$ up to a degree shift by one, i.e.,
modulo degree shift the Morse homologies of $F$ and $f|_{h^{-1}(0)}$ coincide. The proof
of this fact in \cite{frauenfelder} is on homology level and not on chain level. It is still
a desideratum to have a proof of this on chain level in the spirit of an adiabatic limit argument as in \cite{dostoglou-salamon, gaio-salamon, salamon-weber}. While there is a natural one-to-one correspondence between critical points of the two functionals this is 
in general not the case for the gradient flow lines. If $g$ is a Riemannian metric on 
$M$ we consider the product metric $g \oplus g_\mathbb{R}$ on $M \times \mathbb{R}$ 
where $g_\mathbb{R}$ is the standard metric on $\mathbb{R}$. Then gradient flow lines
$(x,\tau) \colon \mathbb{R} \to M \times \mathbb{R}$ are solutions of the ODE
\begin{equation}\label{gr1}
\left\{\begin{array}{c}
\partial_s x(s)+\nabla f(x(s))+\tau(s)\nabla h(x(s))=0 \\
\partial_s \tau(s)+h(x(s))=0
\end{array}\right.
\end{equation}
where $s \in \mathbb{R}$ and $\nabla$ denotes the gradient with respect to the metric
$g$ on $M$. On the other hand the gradient flow equation for the restriction
$f|_{h^{-1}(0)}$ with respect to the restriction of the metric $g$ to $h^{-1}(0)$ reads
\begin{equation}\label{gr2}
\left\{\begin{array}{c}
\partial_s x(s)+\nabla f(x(s))+\tau(s)\nabla h(x(s))=0 \\
h(x(s))=0.
\end{array}\right.
\end{equation}
One can interpolate between the two gradient flow equations as follows. Consider for
$\epsilon>0$ the metric $g \oplus \epsilon g_\mathbb{R}$ on $M \times \mathbb{R}$.
Then gradient flow lines of the Lagrange multiplier functional $F$ with respect to
this metric are solutions of the equation
$$\left\{\begin{array}{c}
\partial_s x(s)+\nabla f(x(s))+\tau(s)\nabla h(x(s))=0 \\
\epsilon\partial_s \tau(s)+h(x(s))=0
\end{array}\right.$$
For $\epsilon=1$ this is (\ref{gr1}) and for $\epsilon=0$ we recover (\ref{gr2}).
\\ \\
In this note we consider the two gradient flow equations in an infinite dimensional set-up.
Rabinowitz action functional is the Lagrange multiplier functional of minus the area functional to the constraint given by vanishing of the mean value of a Hamiltonian. 
Namely consider an exact symplectic manifold $(M,\omega=d\lambda)$ and a smooth function
$H \colon M \to \mathbb{R}$ referred to as the Hamiltonian. If $S^1=\mathbb{R}/\mathbb{Z}$
is the circle we abbreviate by
$$\mathcal{L}=C^\infty(S^1,M)$$
the free loop space of $M$. Rabinowitz action functional is the Lagrange multiplier functional
$$\mathcal{A}^H \colon \mathcal{L} \times \mathbb{R} \to \mathbb{R}, \quad
(u,\tau) \mapsto -\int_{S^1}u^*\lambda+\tau \int_0^1 H(u(t))dt.$$
For the $L^2$-metric on $\mathcal{L}$ obtained by integrating the family of
Riemannian metrics $\omega(\cdot,J_t \cdot)$ on $M$, where $J_t$ for $t \in S^1$
is a time-dependent family of $\omega$-compatible almost complex structures 
gradient flow lines of $\mathcal{A}^H$ are solutions $(u,\tau)
\in C^\infty(\mathbb{R}\times S^1,M) \times C^\infty(\mathbb{R},\mathbb{R})$  of the
equation
\begin{equation}\label{grad1}
\left\{\begin{array}{c}
\partial_s u(s,t)+J_t(u(s,t))\Big(\partial_t u(s,t)-\tau(s)X_H(u(s,t))\Big)=0 \\
\partial_s \tau(s)+\int_0^1 H(u(s,t)dt=0.
\end{array}\right.
\end{equation}
where $X_H$ is the Hamiltonian vector field of $H$ implicitly defined by the condition
$$dH=\omega(\cdot,X_H).$$
The energy of a solution $(u,\tau)$ of (\ref{grad1}) is defined as
$$E(u,\tau)=\int_{-\infty}^\infty \int_0^1 \omega\Big(\partial_s u(s,t),
J_t(u(s,t))\partial_s u(s,t)\Big)dt ds+\int_{-\infty}^\infty (\partial_s \tau)^2 ds.$$
Note that since $J_t$ is $\omega$-compatible the energy is nonnegative but could be
infinite.
We abbreviate by
$$\mathcal{M}_1=\Big\{(u,\tau)\,\,\textrm{solution of (\ref{grad1})},\,\,E(u,\tau)<\infty\Big\}$$
the moduli space of finite energy solutions of (\ref{grad1}).
\\ \\
If we consider the restriction of minus the area functional
$u \mapsto -\int u^*\lambda$ on $\mathcal{L}$ to the constraint given by the vanishing
of the mean value of $H$ on free loops we obtain the gradient flow equation
\begin{equation}\label{grad2}
\left\{\begin{array}{c}
\partial_s v(s,t)+J_t(u(s,t))\Big(\partial_t v(s,t)-\tau(s)X_H(v(s,t))\Big)=0 \\
\int_0^1 H(v(s,t)dt=0.
\end{array}\right.
\end{equation}
In this case $\tau$ is uniquely determined by $v$. For solutions $v$ of (\ref{grad2})
the energy is defined as
$$E(v)=\int_{-\infty}^\infty \int_0^1 \omega\Big(\partial_s v(s,t),
J_t(v(s,t))\partial_s v(s,t)\Big)dt ds$$
and we abbreviate by
$$\mathcal{M}_2=\Big\{v\,\,\textrm{solution of (\ref{grad2})},\,\,E(v)<\infty\Big\}$$
the moduli space of finite energy solutions of (\ref{grad2}).
\\ \\
In this paper we are studying the case where our symplectic manifold
$M=\mathbb{R} \times \Sigma$ is the symplectization of a contact manifold $\Sigma$.
Our Hamiltonian is the map
$$H \colon \mathbb{R} \times \Sigma \to \mathbb{R}, \quad (r,x) \mapsto e^r-1$$
and our family of $\omega$-compatible almost complex structures $J_t$ are
additionally required to be SFT-like. In this case we have a natural map
$$\Psi \colon \mathcal{M}_1 \to  \mathcal{M}_2.$$
Namely given $(u,\tau) \in \mathcal{M}_1$ define
$$\sigma_u \colon \mathbb{R} \to \mathbb{R}$$
by the condition that for $s \in \mathbb{R}$ it holds
$$\int_0^1 H\big(\sigma_u(s)_*u(s,t)\big)dt=0.$$
Here we use the obvious $\mathbb{R}$-action on $\mathbb{R} \times \Sigma$
by translation of the first factor, namely if $r_1,r_2 \in \mathbb{R}$ and
$x \in \Sigma$, then
$$(r_1)_*(r_2,x)=(r_1+r_2,x).$$
Explicitly, $\sigma_u(s)$ can be computed as
$$\sigma_u(s)=-\ln \bigg(\int_0^1 H(u(s,t))dt-1\bigg).$$
With this notion we define
$$\Psi(u,\tau)(s,t):=\sigma_u(s)_*u(s,t), \quad (s,t) \in \mathbb{R} \times S^1.$$
The main result of this paper is the following theorem.
\\ \\
\textbf{Theorem\,A: } \emph{The map $\Psi \colon \mathcal{M}_1 \to \mathcal{M}_2$
is a bijection.}
\\ \\
The proof of Theorem\,A consists of constructing a map
$$\Phi \colon \mathcal{M}_2 \to \mathcal{M}_1$$
and then showing that $\Phi$ is the inverse of $\Psi$. The construction of 
$\Phi$ uses an existence and uniqueness result for a Kazdan-Warner equation. 
The proof of Theorem\,A is carried out in Section~\ref{pro}. 
\\ \\
In the Appendix the main motivation
of the author to study this new version of the gradient flow equation of Rabinowitz
action functional is explained. Namely the restriction of the negative area functional
to the constraint given by the mean value of the Hamiltonian is not only antiinvariant
under time reversal but satisfies as well Chas-Sullivan additivity and therefore
a Floer homology for this action functional should profit from both of these properties on
chain level and not just homology level. 
\\ \\
\emph{Acknowledgements:} The author acknowledges partial support by DFG grant
FR 2637/2-2.

\section{The symplectization of a contact manifold}

In this section we recall the symplectization of a contact manifold to fix notation. 
Assume that $(\Sigma,\lambda)$ is a $2n-1$-dimensional contact manifold, i.e.,
the contact form $\lambda$ is a one-form on $\Sigma$, such that
$$\lambda \wedge d\lambda^{n-1}>0.$$
We denote by $R$ the Reeb vector field of $\lambda$ on $\Sigma$ implicitly defined by
$$d\lambda(R,\cdot)=0, \quad \lambda(R)=1.$$
Abbreviate by 
$$\xi=\ker \lambda$$
the hyperplane plane distribution on $T\Sigma$ referred to as the contact structure. 
The restriction of $d\lambda$ to $\xi$ is symplectic, so that $\xi$ becomes a 
symplectic vector bundle over $\Sigma$ of rank $2n-2$. Choose a $d \lambda$-compatible
almost complex structure $J$ on $\xi$, i.e. $d\lambda(\cdot,J \cdot)$, is a bundle
metric on $\xi$. All these structures on $\Sigma$ have canonical extensions to its
symplectization $\mathbb{R} \times \Sigma$. By abuse of notation we denote these canonical
extensions by the same letters. Namely the one-form $\lambda$ gives rise to a one-form
$\lambda$ on the symplectization by setting
$$\lambda_{r,x}=e^r \lambda_x, \quad (r,x) \in \mathbb{R} \times \Sigma.$$
In particular, if we identify $\Sigma$ with $\{0\} \times \Sigma$ in the symplectization we
recover the contact form by restricting $\lambda$ to $\Sigma$. Note that 
$$\omega_{r,x}=d\lambda_{r,x}=e^r d\lambda_x+e^r dr \wedge \lambda_x$$ 
is a symplectic form on $\mathbb{R} \times \Sigma$. Using the splitting
$$T_{r,x}(\mathbb{R} \times \Sigma)=\mathbb{R} \times T_x \Sigma, \quad 
(r,x) \in \mathbb{R} \times \Sigma$$
we extend the Reeb vector field $R$ and the hyperplane distribution $\xi$ to
$T(\mathbb{R} \times \Sigma)$ trivially on the $\mathbb{R}$-factor. We extend $J$ to
an $\omega$-compatible almost complex structure on $T(\mathbb{R} \times \Sigma)$ by the
requirement that
$$JR=-\partial_r,\quad J\partial_r=R,$$
i.e., $J$ interchanges up to sign the Reeb vector field $R$ and the Liouville vector field
$\partial_r$. Such an $\omega$-compatible almost complex structure on the symplectization
$\mathbb{R} \times \Sigma$ is referred to as an \emph{SFT-like almost complex structure},
since it is invariant under the obvious $\mathbb{R}$-action on $\mathbb{R} \times
\Sigma$ and preserves the symplectic splitting
$$T(\mathbb{R}\times \Sigma)=\xi \oplus \langle R,\partial_r \rangle.$$
We now consider the Hamiltonian
\begin{equation}\label{ham}
H \colon \mathbb{R} \times \Sigma \to \mathbb{R}, \quad (r,x) \mapsto e^r-1.
\end{equation}
Note that 
$$H^{-1}(0)=\Sigma$$
and the Hamiltonian vector field $X_H$ of $H$ implicitly defined by the condition
$dH=\omega(\cdot,X_H)$ just equals the Reeb vector field, i.e.
$$X_H=R.$$

\section{A Kazdan Warner equation}

In this section we discuss existence and uniqueness of solutions of a Kazdan Warner equation
\cite{bieberbach, kazdan-warner}.
In the following we use this solution to construct the map $\Phi$.
\\ \\
We consider a smooth nonegative function $b \colon \mathbb{R} \to [0,\infty)$
with the property that
$$||b||_{L^1}=\int_{-\infty}^\infty b(s)ds<\infty.$$
For the readers convenience we give in this section a proof of the following theorem,
which for experts of Kazdan Warner equations is probably known. 
\begin{thm}\label{kw}
There exists a unique $\rho \in W^{2,2}(\mathbb{R})$ solving the nonlinear 
second order ODE
\begin{equation}\label{ode}
\partial^2_s \rho=1-e^{-\rho}-b.
\end{equation}
Moreover, this unique solution has the property to be nonnegative. 
\end{thm}
Before we can embark on the proof of the theorem we need various preparations.
\begin{lemma}\label{bound}
Suppose that $\rho \in W^{2,2}(\mathbb{R})$ is a solution of (\ref{ode}).
Then for every $s \in \mathbb{R}$ we have
$$0 \leq \rho(s) \leq \max\big\{2\ln 2, 4||b||_{L^1}^2\big\}.$$
\end{lemma}
\textbf{Proof: } We first show that for every $s \in \mathbb{R}$ we have that $\rho(s)$
is nonnegative. To see that we argue by contradiction and assume that there exists
$s$ such that $\rho(s)<0$. Since $\rho \in W^{2,2}(\mathbb{R})$ we have $\lim_{s \to \pm \infty}\rho(s)=0$ and hence there has 
to exist a point $s_0$ at which attains $\rho$ its global negative minimum, in particular
it holds that
$$\rho(s_0)<0, \quad \partial_s^2\rho(s_0) \geq 0.$$
However, since $b$ is nonnegative this contradicts (\ref{ode}).
This shows that $\rho$ is nonnegative. 
\\ \\
To show the upper bound we again argue by contradiction and assume that there exists
$s \in \mathbb{R}$ such that $\rho(s)>\max\{2\ln 2, 4||b||_{L_1}^2\}$. Since $\rho$ asymptotically
converges to zero, it has to attain its global maximum. Therefore there exists
$$\mu>\max\big\{2\ln 2, 4||b||_{L^1}^2\big\}$$
and $s_0 \in \mathbb{R}$ such that
$$\rho(s_0)=\mu, \quad \partial_s \rho(s_0)=0.$$
Since we have already shown that $\rho$ is nonnegative it follows from 
(\ref{ode}) that
$$\partial^2_s \rho \geq -b$$
and therefore for every $s>s_0$ 
\begin{equation}\label{vel}
\partial_s \rho(s)=\int_{s_0}^s\partial^2_s \rho(\sigma)d\sigma \geq 
-\int_{s_0}^s b(\sigma)d\sigma \geq -\int_{-\infty}^\infty b(\sigma)d\sigma=-||b||_{L^1}.
\end{equation}
Define
$$s_1:=\min\big\{s>s_0: \rho(s)=\tfrac{\mu}{2}\big\}$$
the first instant after $s_0$ where $\rho$ attains the value $\tfrac{\mu}{2}$. Note
that since $\rho$ asymptotically converges to zero such an instant has to exist by the
intermediate value theorem. Moreover, since this is the first instant, where 
$\rho$ attains the value $\tfrac{\mu}{2}$ after $s_0$ again, we necessarily have
\begin{equation}\label{nonneg}
\partial_s \rho(s_1) \leq 0.
\end{equation}
Using (\ref{vel}) we estimate
$$\frac{\mu}{2}=-\int_{s_0}^{s_1}\partial_s \rho(s)ds\leq \int_{s_0}^{s_1} ||b||_{L^1}ds
=||b||_{L^1}(s_1-s_0)$$ 
and therefore we have
\begin{equation}\label{zeit}
s_1-s_0 \geq \frac{\mu}{2||b||_{L^1}}>2||b||_{L^1}.
\end{equation}
By definition of $s_1$ we have for every $s \in [s_0,s_1]$
$$\rho(s) \geq \frac{\mu}{2}>\ln 2$$
and therefore in view of (\ref{ode})
$$\partial^2_s \rho(s) \geq \frac{1}{2}-b, \quad s \in [s_0,s_1].$$
This inequality together with (\ref{zeit}) implies
\begin{eqnarray*}
\partial_s \rho(s_1)&=&\int_{s_0}^{s_1} \partial_s^2 \rho(s)ds\\
&\geq& \int_{s_0}^{s_1}\bigg(\frac{1}{2}-b\bigg)ds\\
&=&\frac{s_1-s_0}{2}-\int_{-\infty}^\infty bds\\
&>&||b||_{L^1}-||b||_{L^1}\\
&=&0
\end{eqnarray*}
contradicting (\ref{nonneg}). This proves the lemma. \hfill $\square$
\begin{lemma}
Assume that $\rho \in W^{2,2}(\mathbb{R})$ is a solution of (\ref{ode}), then
\begin{equation}\label{en}
\int_{-\infty}^\infty e^{-\rho}\big(b+(\partial_s \rho)^2\big)ds
+\int_{-\infty}^\infty \big(1-e^{-\rho}\big)^2ds=||b||_{L^1}.
\end{equation}
\end{lemma}
\textbf{Proof: } Using (\ref{ode}) we compute
\begin{equation}\label{kawa}
\partial^2_s(\rho+e^{-\rho})=(1-e^{-\rho})\partial^2_s \rho+e^{-\rho}(\partial_s
\rho)^2=(1-e^{-\rho})^2+e^{-\rho}\big((\partial_s \rho)^2+b\big)-b.
\end{equation}
Since $\rho \in W^{2,2}$ it follows that
$$\lim_{s \to \pm \infty} \partial_s (\rho+e^{-\rho})=\lim_{s \to \pm \infty}
(1-e^{-\rho})\partial_s \rho=0$$
and therefore
$$\int_{-\infty}^\infty \partial^2_s(\rho+e^{-\rho})ds=0.$$
Hence (\ref{en}) follows from integrating (\ref{kawa}). \hfill $\square$
\\ \\
Before stating the next lemma we first point out that since $b$ is smooth and
has finite $L^1$-norm it follows that it has as well finite $L^2$-norm $||b||_{L^2}$. 
\begin{lemma}\label{wbo}
There exists a constant 
$$c=c\big(||b||_{L^1},||b||_{L^2}\big)$$ 
depending only on the
$L^1$-norm and $L^2$-norm of $b$ and which can be chosen to depend continuously
on these two norms such that for every solution $\rho \in W^{2,2}(\mathbb{R})$ of
(\ref{ode}) we have 
$$||\rho||_{W^{2,2}} \leq c.$$
\end{lemma}
\textbf{Proof: } We have
$$
||\rho||_{W^{2,2}}^2=||\partial^2_s \rho||_{L^2}^2+||\partial_s \rho||_{L^2}^2
+||\rho||_{L^2}^2$$
and we estimate all three terms on the righthand side individually. 
In view of (\ref{ode}) and (\ref{en}) it holds that
\begin{equation}\label{w2}
||\partial^2_s \rho||_{L^2} \leq ||1-e^{-\rho}||_{L^2}+||b||_{L^2}
\leq \sqrt{||b||_{L^1}}+||b||_{L^2}.
\end{equation}
By convexity of the exponential function we have for every $\kappa>0$
$$x \leq \frac{\kappa(1-e^{-x})}{1-e^{-\kappa}}, \quad x \in [0,\kappa].$$
In view of Lemma~\ref{bound} there exists therefore a constant $c_1=c_1(||b||_{L^1})$
depending continuously on $||b||_{L^1}$ such that
$$\rho \leq c_1(1-e^{-\rho}).$$
Therefore in view of (\ref{en}) we have
\begin{equation}\label{l2}
||\rho||_{L^2}^2 \leq c_1^2 ||b||_{L^1}.
\end{equation}
Using integration by parts and the Cauchy-Schwarz inequality we finally have
$$||\partial_s \rho||_{L^2}^2 \leq ||\rho||_{L^2}\cdot ||\partial_s^2 \rho||_{L^2}$$
so that a uniform bound on the $L^2$-norm of $\partial_s \rho$ follows from
(\ref{w2}) and (\ref{l2}). This proves the lemma. \hfill $\square$
\\ \\
We next explain that the moduli space of solutions of (\ref{ode}) is always regular.
For that purpose it is useful to interpret the moduli space of solutions as
the zero set of a smooth map between Hilbert spaces. 
Namely we consider
$$\mathcal{F}_b \colon W^{2,2}(\mathbb{R}) \to L^2(\mathbb{R}),\quad
\rho \mapsto \partial^2_s \rho+e^{-\rho}-1+b.$$
Then solutions of (\ref{ode}) correspond to the zeros of the map $\mathcal{F}_b$. 
That the image of $\mathcal{F}_b$ actually lies in $L^2(\mathbb{R})$ is the content of
the following little lemma.
\begin{lemma}
Suppose $\rho \in W^{2,2}(\mathbb{R})$. Then 
$$\partial^2_s \rho+e^{-\rho}-1+b \in L^2(\mathbb{R}).$$
\end{lemma}
\textbf{Proof: } Since $b$ is smooth and in $L^1(\mathbb{R})$ by assumption it lies
in $L^2(\mathbb{R})$. Because $\rho \in W^{2,2}(\mathbb{R})$ it is in particular continuous
and uniformly bounded. Therefore there exists a constant $c=c(\rho)$ not depending
on $s$ such that
$$e^{-\rho(s)} \leq 1+c|\rho(s)|, \quad \forall\,\,s \in \mathbb{R}$$
and hence
$$||e^{-\rho}-1||_{L_2} \leq c||\rho||_{L^2} \leq c||\rho||_{W^{2,2}}.$$
This proves the lemma. \hfill $\square$
\\ \\
For any $\rho \in W^{2,2}(\mathbb{R})$ the differential of $\mathcal{F}_b$ at $\rho$
is the bounded linear operator
$$D_\rho:=d\mathcal{F}_b(\rho) \colon W^{2,2}(\mathbb{R}) \to L^2(\mathbb{R}),
\quad \xi \mapsto \partial_s^2 \xi-e^{-\rho}\xi.$$
From the following proposition it follows that the moduli space of solutions of
(\ref{ode}) is always regular and consists of an isolated set of points. 
\begin{prop}\label{reg}
For any $\rho \in W^{2,2}(\mathbb{R})$ the operator $D_\rho$ is an isomorphism
between $W^{2,2}(\mathbb{R})$ and $L^2(\mathbb{R})$. 
\end{prop}
\textbf{Proof: } We prove the Proposition in four steps.
\\ \\
\textbf{Step\,1:} $D_\rho$ is injective for every $\rho \in W^{2,2}(\mathbb{R})$.
\\ \\
In order to prove Step\,1 we assume
that $\xi$ is in the kernel of $D_\rho$, i.e.
$$D_\rho(\xi)=0.$$
We take the $L^2$-inner product of $D_\rho \xi$ with $\xi$ and obtain via integration
by parts
$$0=\langle D_\rho \xi,\xi\rangle=\int_{-\infty}^\infty (\partial^2_s\xi) \xi ds
-\int_{-\infty}^\infty e^{-\rho}\xi^2 ds=-\int_{-\infty}^\infty (\partial_s\xi)^2 ds
-\int_{-\infty}^\infty e^{-\rho}\xi^2 ds$$
which implies that $\xi=0$ and hence $D_\rho$ is injective.
\\ \\
\textbf{Step\,2: } $D_0 \colon W^{2,2}(\mathbb{R}) \to L^2(\mathbb{R})$ is an isomorphism.
\\ \\
In view of Step\,1 it suffices to show that $D_0$ is surjective. Pick $\eta
\in L^2(\mathbb{R})$. We have to find $\xi \in W^{2,2}(\mathbb{R})$ such that
$$\partial^2_s \xi(s)-\xi(s)=\eta(s),\quad s \in \mathbb{R}.$$
Applying the Fourier-Plancherel transform to this equation (see for instance \cite[p.\,188-189]{rudin}) we obtain for the Fourier-Plancherel transforms $\hat{\xi}$ and
$\hat{\eta}$ of $\xi$ respectively $\eta$ the equation
$$\hat{\eta}(s)=-s^2 \hat{\xi}(s)-\hat{\xi}(s), \quad s \in \mathbb{R}.$$
Abbreviating
$$\phi \colon \mathbb{R} \to \mathbb{R}, \quad s \mapsto -\frac{1}{1+s^2}$$
we can rewrite this as
$$\hat{\xi}=\phi \cdot \hat{\eta}.$$
Applying the Fourier-Plancherel transform once more to this identity, we finally
define $\xi$ by 
$$\xi(s):=\widehat{\phi \cdot \hat{\eta}}(-s), \quad s \in \mathbb{R}.$$
Then
$$D_0\xi=\eta$$
and Step\,2 is proved. 
\\ \\
\textbf{Step\,3: } $D_\rho$ is a Fredholm operator of index zero for every 
$\rho \in W^{2,2}(\mathbb{R})$. 
\\ \\
It follows from Step\,2 that $D_0$ is a Fredholm operator of index zero. Since the propery
of being Fredholm is open, there exists $\epsilon>0$ such that each bounded operator
$D \colon W^{2,2}(\mathbb{R}) \to L^2(\mathbb{R})$ satisfying
$$||D-D_0||<\epsilon$$
is Fredholm of index zero. Here the norm denotes the operator norm. Suppose that
$\rho \in W^{2,2}(\mathbb{R})$ choose a smooth cut-off function
$\beta \colon \mathbb{R} \to [0,1]$ with the property that 
$$\beta(s)=1, \quad |s| \geq T$$
for some $T>0$ satisfying
$$||1-e^{-\beta \rho}||_{L^\infty} <\epsilon.$$
We have the difference of operators
$$D_{\beta \rho}-D_0 \colon W^{2,2}(\mathbb{R}) \to L^2(\mathbb{R}), \quad
\xi \mapsto (1-e^{-\beta \rho})\xi.$$
Note that
$$||(D_{\beta \rho}-D_0)\xi||_{L^2}\leq ||1-e^{-\beta \rho}||_{L^\infty} 
||\xi||_{L^2} \leq \epsilon ||\xi||_{W^{2,2}}$$ 
so that we have
$$||D_{\beta \rho}-D_0|| \leq \epsilon$$
and therefore $D_{\beta \rho}$ is a Fredholm operator of index zero. We further have the
difference of operators
$$D_\rho-D_{\beta \rho} \colon W^{2,2}(\mathbb{R}) \to L^2(\mathbb{R}), \quad
\xi \mapsto (e^{-\beta \rho}-e^{-\rho})\xi.$$
Since $\beta$ equals one outside a compact subset of $\mathbb{R}$ the
continuous function $e^{-\beta \rho}-e^{-\rho}$ has compact support and therefore
the operator $D_\rho-D_{\beta \rho}$ is compact. Since adding a compact operator to
Fredholm operator gives still a Fredholm operator of the same index we conclude that
$D_\rho$ is a Fredholm operator of index zero as well. 
\\ \\
\textbf{Step\,4: } We prove the Proposition.
\\ \\
By Step\,3 we have
$$0=\mathrm{ind} D_\rho=\dim \ker D_\rho-\dim \mathrm{coker} D_\rho.$$
From Step\,1 we know that $D_\rho$ is injective and therefore
$$\dim \ker D_\rho=0.$$
Therefore we have
$$\dim \mathrm{coker} D_\rho=0$$
so that $D_\rho$ is surjective as well and therefore an isomorphism. This proves
the Proposition. \hfill $\square$
\\ \\
We are finally in position to prove the main result of this section.
\\ \\
\textbf{Proof of Theorem~\ref{kw}: } We first discuss uniqueness. For that purpose
let us assume that $\rho_1, \rho_2 \in W^{2,2}(\mathbb{R})$ both solve (\ref{ode}),\,i.e.,
$$\partial_s^2 \rho_1+e^{-\rho_1}-1=b=\partial_s^2 \rho_2+e^{-\rho_2}-1.$$
Their difference solves the second order ODE
\begin{equation}\label{dode}
\partial^2_s(\rho_2-\rho_1)=e^{-\rho_2}\big(e^{\rho_2-\rho_1}-1\big).
\end{equation}
Since $\rho_2-\rho_1 \in W^{2,2}(\mathbb{R})$ it converges asymptotically to zero. 
We want to show that it is identically zero. To see that we argue by contradiction.
Otherwise if
$\rho_2-\rho_1$ were not constant zero there exists a local maximum of $\rho_2-\rho_1$ at which
$$(\rho_2-\rho_1)(s_0)>0, \quad \partial_s^2 (\rho_2-\rho_1)(s_0) \leq 0$$
or a local minimum at which
$$(\rho_2-\rho_1)(s_0)<0, \quad \partial_s^2 (\rho_2-\rho_1)(s_0) \geq 0.$$
Both of these contradict (\ref{dode}) and therefore $\rho_2-\rho_1=0$, i.e.,
$\rho_2=\rho_1$, and uniqueness follows.
\\ \\
It remains to prove existence. For that purpose we introduce the following subset
of the closed interval
$$E=\big\{r \in [0,1]: \mathcal{F}_{rb}^{-1}(0) \neq \emptyset\big\},$$
i.e., the set of all $r \in [0,1]$ for which there exists a solution 
$\rho \in W^{2,2}(\mathbb{R})$ of the ODE
\begin{equation}\label{rode}
\partial^2_s \rho=1-e^{-\rho}-rb.
\end{equation}
We first observe that
$0 \in E$, since indeed $\rho=0$ is a solution of (\ref{rode}) for $r=0$. We next
discuss that $E$ is a closed subset of the interval $[0,1]$. For that purpose
suppose that $r_\infty \in [0,1]$ and there exists a sequence $r_\nu \in E$ with
$\nu \in \mathbb{N}$ such that 
$$\lim_{\nu \to \infty} r_\nu=r_\infty.$$
Since $r_\nu \in E$, there exists $\rho_\nu \in W^{2,2}(\mathbb{R})$ such that
$\rho_\nu$ is a solution of (\ref{rode}) for $r=r_\nu$. By Lemma~\ref{wbo} the
$W^{2,2}(\mathbb{R})$ norm of the sequence $\rho_\nu$ is uniformly bounded.
Therefore by the Theorem of Banach Alaoglu there exists $\rho \in W^{2,2}(\mathbb{R})$
such that $\rho_\nu$ converges weakly to $\rho$. The map $\rho$ is then a solution of
(\ref{rode}) for $r=r_\infty$. In particular, $r_\infty \in E$ which shows that
$E$ is closed. We finally note that $E$ is open in view of Proposition~\ref{reg}
and the Implicit Function theorem. We have checked that $E$ is a nonempty, open and
closed subset of the interval $[0,1]$ and since the interval is connected it follows that
$$E=[0,1].$$
In particular, $1 \in E$ and this means that there exists a solution 
$\rho \in W^{2,2}(\mathbb{R})$ of (\ref{rode}) for $r=1$. This proves existence. That the solution is
nonnegative follows from Lemma~\ref{bound}. The theorem is proven. \hfill $\square$
\section{Construction of the map $\Phi$}

In this section we assume that $M=\mathbb{R} \times \Sigma$ is the symplectization
of a contact manifold $\Sigma$, the Hamiltonian $H \colon M \to \mathbb{R}$ is
given by (\ref{ham}) and $J_t$ is a smooth family of SFT-like almost complex structures
on $M$. The following lemma tells us that for solutions $v$ of (\ref{grad2}) the Lagrange multiplier is just given by the area of $v$. In particular, since $v$ is a gradient flow
line of minus the area functional the Lagrange multiplier is a monotone increasing function. 
\begin{lemma}\label{area}
Suppose that $v$ is a solution of (\ref{grad2}). Then the Lagrange multiplier satisfies
$$\tau(s)=\int_0^1 v_s^* \lambda$$
for every $s \in \mathbb{R}$, where $v_s=v(s,\cdot) \colon S^1 \to M$.
\end{lemma}
\textbf{Proof: } We compute using (\ref{grad2})
\begin{eqnarray*}
0&=&\int_0^1 dH(v)\partial_s v dt\\
&=&-\int_0^1 dH(v) J_t(v)\partial_t v dt+\tau \int_0^1 dH(v) J_t(v)X_H(v)dt\\
&=&-\int_0^1 \omega(J_t(v)\partial_t v,X_H)dt-\tau \int_0^1 dH(v)\partial_r dt\\
&=&-\int_0^1 \omega(\partial_t v,\partial_r)dt-\tau \int_0^1\big(H(v)+1\big)dt\\
&=&\int_{S^1} v^* \lambda-\tau.
\end{eqnarray*}
This proves the lemma. \hfill $\square$
\\ \\
Suppose now that $v \in \mathcal{M}_2$. We define a smooth function
$$b_v \colon \mathbb{R} \to \mathbb{R}, \quad s \mapsto \partial_s \int_0^1 v_s^* \lambda.$$
Since $v$ is a gradient flow line of minus the area functional we have
$$b_v(s) \geq 0$$
for every $s \in \mathbb{R}$. Moreover, the $L^1$-norm of $b_v$ is given by the the energy
of $v$,\,i.e.,
$$||b_v||_{L^1}=E(v).$$
By Theorem~\ref{kw} there exists a unique solution $\rho_v \in W^{2,2}(\mathbb{R})$ of the ODE
\begin{equation}\label{vode}
\partial^2_s \rho_v=1-e^{-\rho_v}-b_v.
\end{equation}
Define 
$$u_v:=(-\rho_v)_* v \colon \mathbb{R} \times S^1 \to M.$$
We claim that the tuple
$$w_v:=\bigg(u_v,\int_0^1 v^*\lambda+\partial_s \rho_v\bigg)$$
is a solution of the gradient flow equation (\ref{grad1}). We compute
using Lemma~\ref{area} and taking advantage that SFT-like almost complex structures
are invariant under the $\mathbb{R}$-action on $M=\mathbb{R} \times \Sigma$
\begin{eqnarray*}
\partial_s u_v&=&\partial_s (-\rho_v)_* v\\
&=&d(-\rho_v)_* \partial_s v-(\partial_s \rho_v)\partial_r\\
&=&-d(-\rho_v)_* J(v)\Bigg(\partial_t v-\bigg(\int_0^1 v^*\lambda\bigg) X_H(v)\Bigg)+
(\partial_s \rho_v)J(u_v) X_H(u_v)\\
&=&-J(u_v)\Bigg(\partial_t u_v-\bigg(\int_0^1 v^*\lambda\bigg)X_H(u_v)\Bigg)
+(\partial_s \rho_v)J(u_v) X_H(u_v)\\
&=&-J(u_v)\Bigg(\partial_t u_v-\bigg(\int_0^1 v^*\lambda+\partial_s \rho_v
\bigg)X_H(u_v)\Bigg).
\end{eqnarray*}
This proves the first equation of the gradient flow equation (\ref{grad1}). To check
the second equation in (\ref{grad1}) we first note that for $(r,x) \in \mathbb{R}\times
\Sigma$ and $\rho \in \mathbb{R}$ we have
\begin{eqnarray*}
H\big((-\rho_*)(r,x)\big)&=&H(r-\rho,x)\\
&=&e^{r-\rho}-1\\
&=&e^{-\rho}\big(e^r-e^\rho\big)\\
&=&e^{-\rho}\big(H(r,x)+1-e^\rho\big)\\
&=& e^{-\rho}H(r,x)+e^{-\rho}-1
\end{eqnarray*}
and therefore we compute using (\ref{vode}) and the fact that the mean value of
$H(v)$ vanishes according to (\ref{grad2})
\begin{eqnarray*}
\partial_s \bigg(\int_0^1 v^*\lambda+\partial_s \rho_v\bigg)&=&b_v+\partial^2_s \rho_v\\
&=&1-e^{-\rho_v}\\
&=&-e^{-\rho_v}\int_0^1 H(v)dt+1-e^{-\rho_v}\\
&=&-\int_0^1 H\big((-\rho_v)_* v\big)dt\\
&=&-\int_0^1 H(u_v)dt.
\end{eqnarray*}
Therefore the second equation in (\ref{grad1}) holds true as well and we define
$$\Phi \colon \mathcal{M}_2 \to \mathcal{M}_1$$
for $v \in \mathcal{M}_2$ by
$$\Phi(v)=w_v.$$

\section{Proof of the main result}\label{pro}

Theorem\,A from the Introduction now follows from the following result.
\begin{thm}
The map $\Phi$ is inverse to the map $\Psi$.
\end{thm}
\textbf{Proof: } We prove the theorem in two steps.
\\ \\
\textbf{Step\,1:} We have $\Psi \circ \Phi=\mathrm{id} \colon \mathcal{M}_2 \to \mathcal{M}_2$, i.e.,
$\Phi$ is right inverse to $\Psi$. 
\\ \\
Suppose that $v \in \mathcal{M}_2$. By construction both maps $\Psi$ and $\Phi$
act via the $\mathbb{R}$-action on $\mathbb{R} \times \Sigma$. Therefore there
exists a smooth function $\chi \colon \mathbb{R} \to \mathbb{R}$ such that
$$\Psi \circ \Phi(v)=\chi_* v.$$
Since both $v$ and $\chi_* v$ belong to the moduli space $\mathcal{M}_2$ we have
$$0=\int_0^1 H(\chi_* v)dt=e^\chi \int_0^1 H(v)dt+e^\chi-1=e^\chi-1$$
and therefore
$$\chi=0.$$
This proves that
$$\Psi \circ \Phi(v)=v$$
and hence $\Phi$ is right inverse to $\Psi$. 
\\ \\
\textbf{Step\,2:} We have $\Phi \circ \Psi=\mathrm{id} \colon \mathcal{M}_1 
\to \mathcal{M}_1$, i.e., $\Phi$ is left inverse to $\Psi$.
\\ \\
Suppose that $(u,\tau) \in \mathcal{M}_1$. We abbreviate
$$v:=\Psi(u,\tau)=(\sigma_u)_* u$$
where
$$\sigma_u=-\ln \bigg(\int_0^1 H(u)dt-1\bigg).$$
As solution of $\mathcal{M}_2$ the Lagrange multiplier for $v$ is completely
determined by $v$ according to Lemma~\ref{area}. We can alternatively express it
as well with the help of $\tau$ and $\sigma_u$. Indeed, in view of
$v=(\sigma_u)_* u$ and the fact that $(u,\tau)$ is a solution of
(\ref{grad1}) we obtain the formula
$$\partial_s v+J(v)\Big(\partial_t v-(\tau-\partial_s \sigma_u)X_H(v)\Big)=0.$$
Therefore in view of Lemma~\ref{area} we have
$$\int_0^1 v^*\lambda=\tau-\partial_s \sigma_u.$$
For the unique solution $\rho_v \in W^{2,2}(\mathbb{R})$ of equation (\ref{vode})
we have
$$\Phi(v)=\bigg((-\rho_v)_*v,
\int_0^1 v^*\lambda+\partial_s \rho_v\bigg)=
\bigg((\sigma_u-\rho_v)_*u, \tau-\partial_s \sigma_u+\partial_s \rho_v\bigg).$$
We abbreviate
$$\chi:=\sigma_u-\rho_v \colon \mathbb{R} \to \mathbb{R}.$$
With this notion it holds that
\begin{equation}\label{leftinv}
\big((\chi)_* u, \tau-\partial_s \chi\big)=\Phi \circ \Psi\big(u,\tau\big).
\end{equation}
Since both $(u,\tau)$ and $\Phi \circ \Psi(u,\tau)$ are solutions of (\ref{grad1}) we compute
\begin{eqnarray}\label{kw2}
\partial^2_s \chi&=&\partial_s \tau -\partial_s (\tau-\partial_s \chi)\\ \nonumber
&=&-\int_0^1 H(u)dt+\int_0^1 H(\chi_* u)dt\\ \nonumber
&=&-\int_0^1 H(u)dt+e^\chi \int_0^1 H(u)dt+e^\chi-1\\ \nonumber
&=&\big(e^\chi-1\big)\bigg(\int_0^1 H(u)dt+1\bigg).
\end{eqnarray}
We claim that the only solution of this problem is
\begin{equation}\label{vanish}
\chi=0.
\end{equation}
To see that we first note that since $H$ takes values in $(-1,\infty)$ we have
\begin{equation}\label{pos}
\int_0^1 H(u)dt+1>0.
\end{equation}
Since both $(u,\tau)$ and $\Phi \circ \Psi(u,\tau)$ have finite energy we must have
$$\lim_{s \to \pm \infty}\chi(s)=0.$$
Hence if $\chi$ did not vanish identically it would attain
a positive local maximum or a negative local minimum, i.e., there would exist $s_0 \in \mathbb{R}$ such that
$$\chi(s_0)>0, \quad \partial_s^2 \chi(s_0) \leq 0$$
or
$$\chi(s_0)<0, \quad \partial_s^2 \chi(s_0) \geq 0,$$
both contradicting (\ref{kw2}) in view of (\ref{pos}). This proves (\ref{vanish}). Plugging
this equation into (\ref{leftinv}) we obtaint
$$\Phi \circ \Psi(u,\tau)=(u,\tau).$$
This proves Step\,2 and hence the Theorem follows. \hfill $\square$

\appendix

\section{Symmetries and Chas-Sullivan additivity}\label{cha}

In this appendix we give some motivation for exploring the restriction of the area
functional to the constraint given by the mean value of the Hamiltonian. We explain that
the restriction has the same transformation behaviour under the symmetries of the
free loop space as Rabinowitz action functional but additionally satisfies Chas-Sullivan additivity which the Rabinowitz action functional does not. 
\\ \\
If $M$ is a manifold, we have an $S^1$-action on the free loop space $\mathcal{L}=C^\infty(S^1,M)$ by reparametrization, namely if $u \in \mathcal{L}$ and $r \in S^1$
$$r_*u(t)=u(t+r), \quad t \in S^1.$$
Moreover, we have an involution 
$$I \colon \mathcal{L} \to \mathcal{L}, \quad u \mapsto u^-$$
where $u^-$ is the loop traversed backward
$$u^-(t)=u(-t), \quad t \in S^1.$$
Combining the $S^1$-action with the involution $I$ we obtain an action of
$$O(2)=S^1 \ltimes \mathbb{Z}/2\mathbb{Z}.$$
We further have an action of the monoid $\mathbb{N}$ on $\mathcal{L}$ by iteration.
Namely if $n \in \mathbb{N}$ and $u \in \mathcal{L}$, we set
$$n_* u(t)=u(nt), \quad t \in S^1.$$
Again the $\mathbb{N}$-action combines with the $O(2)$-action to an action of their
semidirect product $O(2) \ltimes \mathbb{N}$. We extend these actions to
$\mathcal{L} \times \mathbb{R}$ as follows. If $(u,\tau) \in \mathcal{L} \times \mathbb{R}$
\begin{eqnarray*}
r_*(u,\tau)&=&(r_*u,\tau), \quad r \in S^1\\
I(u,\tau)&=&(Iu,-\tau)\\
n_*(u,\tau)&=&(n_*u,n\tau), \quad n \in \mathbb{N}.
\end{eqnarray*}
Suppose now that $(M,\omega=d\lambda)$ is an exact symplectic manifold and
$H \colon M \to \mathbb{R}$ is a smooth function. Then Rabinowitz action
functional
$$\mathcal{A}^H \colon \mathcal{L} \times \mathbb{R}, \quad
(u,\tau) \mapsto -\int u^*\lambda+\tau\int_0^1 H(u)dt$$
has the following transformation behaviour under these symmetries. If
$(u,\tau) \in \mathcal{L} \times \mathbb{R}$, then
\begin{eqnarray*}
\mathcal{A}^H\big(r_*(u,\tau)\big)&=&\mathcal{A}^H\big(u,\tau\big),\quad r\in S^1\\
\mathcal{A}^H\big(I(u,\tau)\big)&=&-\mathcal{A}^H\big(u,\tau\big)\\
\mathcal{A}^H\big(n_*(u,\tau)\big)&=&n\mathcal{A}^H\big(u,\tau\big), \quad n \in \mathbb{N}.
\end{eqnarray*}
The restriction of minus the area functional to the constraint given by the mean value
of $H$ satisfies the same transformation behaviour under these symmetries. Namely
abbreviate
$$\mathcal{L}_H=\mathcal{H}^{-1}(0)$$
where
$$\mathcal{H} \colon \mathcal{L} \to \mathbb{R}, \quad
u \mapsto \int_0^1 H(u)dt$$ 
is the mean value of $H$ along a loop, we set
$$\mathfrak{a}^H \colon \mathcal{L}_H \to \mathbb{R}, \quad u \mapsto -\int u^*\lambda.$$
Note that all the symmetries we discussed on the free loop space $\mathcal{L}$ keep
the subspace $\mathcal{L}_H$ invariant and $\mathfrak{a}^H$ transforms as
\begin{eqnarray*}
\mathfrak{a}^H(r_*u)&=&\mathfrak{a}^H(u),\quad r\in S^1\\
\mathfrak{a}^H(Iu)\big)&=&-\mathfrak{a}^H(u)\\
\mathfrak{a}^H(n_*u)&=&n\mathfrak{a}^H(u), \quad n \in \mathbb{N}.
\end{eqnarray*}
A distinguishing feature of $\mathfrak{a}^H$ compared to $\mathcal{A}^H$ however, is
its behaviour under concatenation of loops. Suppose that $u,v \in \mathcal{L}$
have the same starting and endpoint, i.e., $u(0)=v(0)$. In this case we define there
concatenation
$$u \# v(t)=\left\{\begin{array}{cc}
u(2t) & 0\leq t \leq \frac{1}{2}\\
v(2t-1) & \frac{1}{2} \leq t \leq 1.
\end{array}\right.$$
Since the $u$ and $v$ have the same starting and ending point their concatenation is
continuous at $t=\tfrac{1}{2}$. On the other hand, in general it is not smooth at
$t=\tfrac{1}{2}$. However, we can interpret 
$u \# v$ as an element in $\mathcal{L}^{1,2}=W^{1,2}(S^1,M)$ the Hilbert manifold
of $W^{1,2}$-loops. Note that the functionals $\mathfrak{a}^H$ and $\mathcal{A}^H$
canonically extend to $W^{1,2}$-loops. The concatenation product 
gives rise to string topology on the free loop space as discovered by Chas and Sullivan
\cite{chas-sullivan}. 
\\ \\
Denote by $\mathcal{L}^{1,2}_H$ the $W^{1,2}$-loops for which the mean value of
the Hamiltonian $H$ vanishes. If $u,v \in \mathcal{L}^{1,2}_H$ are two loops 
with common starting and endpoint their concatenation $u\#v$ still lies in
$\mathcal{L}^{1,2}_H$. The functional $\mathfrak{a}^H$ is \emph{Chas-Sullivan additive}
with respect to concatenation in the following sense
$$\mathfrak{a}^H(u \# v)=\mathfrak{a}^H(u)+\mathfrak{a}^H(v).$$
For the Rabinowitz action functional $\mathcal{A}^H$ this is not true in general. In fact the
question for Rabinowitz action functional is a bit more subtle since one needs
to specify how the Lagrange multiplier transforms under concatenation. Suppose that
$(u,\tau), (v,\sigma) \in \mathcal{L}^{1,2} \times \mathbb{R}$ satisfy
$u(0)=v(0)$. We want to define
$$\rho=\rho(u,v,\tau,\sigma) \in \mathbb{R}$$
such that Chas-Sullivan additivity holds, i.e.,
$$\mathcal{A}^H(u\#v,\rho)=\mathcal{A}^H(u,\tau)+\mathcal{A}^H(v,\sigma).$$
Since the area functional is Chas-Sullivan additive this leads to the requirement
$$\tau \int_0^1 H(u)dt+\sigma \int_0^1 H(v)dt=\rho \int_0^1 H(u\#v)dt
=\rho \bigg(\int_0^1 H(u)dt+\int_0^1 H(v)dt\bigg)$$
and therefore we need to define $\rho$ by
$$\rho=\frac{\tau \int_0^1 H(u)dt+\sigma \int_0^1 H(v)dt}{\int_0^1 H(u)dt+\int_0^1 H(v)dt}.$$
However, this is ill-defined when $\int_0^1 H(u)dt+\int_0^1 H(v)dt=0$.
\\ \\
To the authors knowledge so far nobody directly defined product on Rabinowitz
Floer homology, but products can be defined on homologies isomorphic to Rabinowitz Floer
homology. So was it proved by Abbondandolo and Merry that Rabinowitz Floer homology
is isomorphic to Floer homology on the time-energy extended phase space on which
products can be defined \cite{abbondandolo-merry}. Alternatively following Cieliebak and Oancea one can define products on V-shaped symplectic homology \cite{cieliebak-oancea}. 
That V-shaped symplectic homology
is isomorphic to Rabinowitz Floer homology was proved in \cite{cieliebak-frauenfelder-oancea}. 
\\ \\
That critical points of Rabinowitz action functional are not just periodic orbit
going forward in time but as well backwards distinguishes it from Symplectic homology
\cite{cieliebak-floer-hofer, viterbo} or Symplectic Field theory \cite{eliashberg-givental-hofer} and has interesting connection to Poincar\'e duality \cite{cieliebak-hingston-oancea}
and Tate homology \cite{albers-cieliebak-frauenfelder}.
A distinguishing feature of the action functional $\mathfrak{a}^H$ is that it
satisfies both Chas-Sullivan additivity and antiinvariance under time-reversal and therefore
it should be possible to use this functional to define some Tate version of a Fukaya category
having as objects a class of Legendrians of a contact manifold. Thinking of such a Fukaya category as a kind of mathematical way for making sense of path integrals the Tate property becomes very reminiscent of the Feynman-Stueckelberg interpretation of a positron being
an electron going backwards in time \cite{feynman, stueckelberg}.

\end{document}